\documentclass[12pt]{article} 

\usepackage[left= 0.8in, textwidth=6.9in, textheight=9.0in, top=1.0in]{geometry}

\usepackage{graphicx}
\usepackage{amsmath,amssymb,verbatim}
\usepackage{color}                                    

\newcommand{\EQ}{\begin{equation}\begin{array}{lllllllll}}
\newcommand{\EE}{\end{array}\end{equation}}
\newcommand{\MT}{\left[ \begin{array}{ccccccccc}}
\newcommand{\EM}{\end{array}\right]}
\newcommand{\bean}{\begin{equation}\begin{array}{rcllllllll}}
\newcommand{\eean}{\end{array}\end{equation}}
\newcommand{\bea}{$$\begin{array}{rcllllllll}}
\newcommand{\eea}{\end{array}$$}

\def\Fr{\ds \frac}
\def\ds{\displaystyle}
\def\xsig{x^\sigma}
\def\xb{x^b}
\def\yb{y^b}
\def\xa{x^a}
\def\Pb{P^b}
\def\Pa{P^a}
\def\Xa{X^a}
\def\Xb{X^b}
\def\Yb{Y^b}
\def\IP{{\cal I}}
\def\ISP{{\cal I}^{\sigma}}
\def\IBP{{\cal I}^{b}}
\def\spxb{\underline{x}^b}

\def\spPb{\underline{P}^b}
\def\spPa{\underline{P}^a}
\def\spXa{\underline{X}^a}

\title{\LARGE \bf
Sparsity-Based Kalman Filters for Data Assimilation\thanks{This work was supported in part by U.S. Naval Research Laboratory - Monterey, CA}
}

\author{Wei Kang\thanks{Wei Kang is with Faculty of Applied Mathematics, Naval Postgraduate School, 
        Monterey, CA, USA,
        {\tt\small wkang@nps.edu}}  and 
         Liang Xu\thanks{Liang Xu is with Naval Research Laboratory, Monterey, CA, USA,
         {\tt\small liang.xu@nrlmry.navy.mil}}%
}

\begin{document}

\maketitle
\thispagestyle{empty}
\pagestyle{empty}

\begin{abstract}
Several variations of the Kalman filter algorithm, such as the extended Kalman filter (EKF) and the unscented Kalman filter (UKF), are widely used in science and engineering applications. In this paper, we introduce two algorithms of sparsity-based Kalman filters, namely the sparse UKF and the  progressive EKF. The filters are designed specifically for problems with very high dimensions. Different from various types of ensemble Kalman filters (EnKFs) in which the error covariance is approximated using a set of dense ensemble vectors, the algorithms developed in this paper are based on sparse matrix approximations of error covariance. The new algorithms enjoy several advantages. The error covariance has full rank without being limited by a set of ensembles. In addition to the estimated states, the algorithms provide updated error covariance for the next assimilation cycle. The sparsity of error covariance significantly reduces the required memory size for the numerical computation. In addition, the granularity of the sparse error covariance can be adjusted to optimize the parallelization of the algorithms. 
\end{abstract}

\section{INTRODUCTION}
For dynamical systems with numerical models, data assimilation is an estimation process of combining observational data with a numerical model to obtain an estimate of the system's state. Data assimilation is essential to numerical weather prediction (NWP). The state estimate can be considered as an interpolation of the sparse observational data; and it is used as the initial condition for the numerical forecast process. If the dimension is relatively low and the data set is small, various linear and nonlinear estimators can be found in the literature that have optimal or suboptimal performances. However, to assimilate big data sets with models that have high dimensions, such as those in operational NWP systems with tens of millions of variables, achieving reliable state estimates and error probability distributions is a challenging problem that have been studied for decades with a huge literature. 

There are two categories of methods widely used in NWP, namely variational methods and the ensemble Kalman filter (EnKF) \cite{xu,houtekamer}. The former is based on a weighted least-square optimization, such as the four dimensional variational data assimilation (4D-Var) in a fixed time window or the three dimensional version (3D-Var) that excludes the time variable. The EnKF algorithm is based on the Kalman filter except that the error covariance is approximated using a set of state ensembles. 4D-Var methods are used in operational NWP systems by many meteorological centers. While it serves as an effective method of data assimilation, 4D-Var algorithms have difficulty to explicitly track the evolution of error covariance within its estimation process due to outrageous computational costs and input/output (I/O) loads required by matrices of extremely high dimensions. EnKF, on the other hand, updates information about the error covariance in the form of ensembles. However, the number of ensemble states is significantly smaller than the number of state variables. As a result, the rank deficiency of error covariance tends to deteriorate the integrity of the estimation process unless remedies to the algorithm, such as localization and covariance inflation, are applied.  

Different types of Kalman filters have been developed and widely used in science and engineering applications, such as the EnKF, the extended Kalman filter (EKF) and the unscented Kalman filter (UKF). In this paper, we introduce two algorithms of sparsity-based Kalman fillters, namely the sparse UKF and the  progressive EKF. The goal of the work is to explore innovative ideas that take the advantage of the (almost) sparsity structure of matrices so that analysis and error covariance can be updated effectively and efficiently without the drawback of rank deficiency. The granularity of subproblems for the purpose of algorithm parallelization is also emphasized in the method. The filters are developed specifically for problems with very high dimensions. Different from EnKFs in which the error covariance is represented using a set of dense ensemble vectors, the new algorithms in this paper are based on a sparse but full rank matrix as an approximation of the error covariance. This is made possible because of two assumptions: (a) the error covariance is approximately a sparse matrix; (b) the system model is component based, i.e. the state vectors are divided into components that can be computed independently in parallel.

\section{Sparse UKF}
Consider  a dynamical system model in which the state variable is $x(t)$, where $t=1, 2, 3,\cdots$ represents time steps. The value of observation at $t=k$ is denoted by $y(k)$. The system model is defined as follows,
\EQ
\label{eq_mdl}
 \begin{aligned}
x (k)&={\cal M}(x(k-1))+\eta_{k-1}, && x(k), \eta_{k-1} \in \mathbb{R}^n, \\
y(k)&={\cal H}(x(k))+\delta_k, && y_k, \delta_k \in \mathbb{R}^m, 
\end{aligned}
\EE
where $ \eta_{k-1}$ is a random variable representing the model error. Its covariance is $Q$. The observational error, $ \delta_k$,  has a covariance $R$. In data assimilation, the goal is to estimate the value of $x(k)$ given the observations $y(1), y(2), \cdots, y(k)$ and the model (\ref{eq_mdl}).  If (\ref{eq_mdl}) is linear and if all random variables are Gaussian, then the Kalman filter is an optimal state estimator. For nonlinear systems with non-Gaussian random errors, various types of Kalman filters exist in the literature with successful applications in science and engineering. If a system has a very high dimension, the conventional form of Kalman filter based on a dense error covariance is not applicable. In this section, we introduce an algorithm that is a variation of UKF for problems with approximately sparse error covariances. 

In a sparse matrix/vector, most entries are zeros. In this paper, we use an underbar to emphasize that a vector or matrix is sparse, for instance $\underline{P}$ and $\underline{x}$. A sparse vector is associated with an index set, denoted by ${\cal I}$, consisting of the indices of nonzero entries. For a sparse matrix, columns may have different numbers of nonzero entries. The largest such number is denoted by $N_{sp}$. In sparsity-based algorithms, a full model evaluation is not always necessary. Using a {\it component-based} model can significantly reduce the computational load. In the notation, a component-based model has three inputs: the sparse state variable, its index, and the index of the output state. More specifically, 
\EQ
\label{ad_model}
\underline{x}(k)={\cal M}(\underline{x}(k-1); {\cal I}_1; {\cal I}_2),\\
\EE
where ${\cal I}_1$ is the index set of the sparse vector $\underline{x}(k-1)$ and ${\cal I}_2$ is the index set of $\underline{x}(k)$. The model evaluates only the entries with indices in ${\cal I}_2$, setting all other entries as zeros. Note that $\underline{x}(k)$ is different from the full state variable $x(k)$ because the later is, in general, a dense vector with mostly nonzero entries. Therefore, it is important to specify the index set ${\cal I}_2$ of the sparse vector $\underline{x}(k)$ to be evaluated using a component-based model. For simplicity, we often omit ${\cal I}_1$ in the notation, i.e.
\EQ
\underline{x}(k)={\cal M}(\underline{x}(k-1); {\cal I}),
\EE
where ${\cal I}$ is the same as ${\cal I}_2$ in (\ref{ad_model}). Algebraic operations between sparse vectors, such as addition and dot product, are defined in the same way as dense vectors. Thus, we may conduct operations between sparse vectors and dense vectors, such as adding a sparse vector to a dense vector $\underline{x}+y$ as long as both vectors have the same dimension. A new operation, called merging, between a sparse vector and a dense vector is defined as follows,
\EQ
z= \underline{x}\triangleright y, &
\left\{ \begin{array}{lll}
i\mbox{th component of }z= i\mbox{th component of } x, & \mbox{if } i\in {\cal I}.\\
i\mbox{th component of }z= i\mbox{th component of } y , & \mbox{if } i\not\in {\cal I}.
\end{array}\right.
\EE
If an operation has an underbar, it means that the evaluation is carried out only at a given index set. For instance, given two sparse matrices $\underbar{A}$ and $\underbar {B}$, then $\underbar{A}* \underbar {B}$ is a different matrix from $\underline{\underline{A}* \underline {B}}$. The former is the conventional matrix multiplication between two sparse matrices; the later is a matrix multiplication in which the entries in a given index set are evaluated and all other entries are set to be zeros. Other operations, such as $\underline{\sqrt{\underline{P}}}$, are defined similarly. A summary of notations is summarized in the following table.

\begin{table}[h]
\begin{center}
\begin{tabular}{|c|c|c|c|c|c|}
\hline
{\small Notation}& {\small Definition}& {\small Notation}&{\small Definition} \\
\hline
$x$& state variable&$y$&{\small observation variable} \\
\hline
 ${\cal M}$&{\small model function}&${\cal H}$&{\small observation operator}\\
 \hline
$n$& {\small state space dimension }& $t=1,2...$&{\small (discrete) time variable}\\
\hline
$\xsig_i$&{\small $\sigma$-point at $t=k-1$}&&\\
\hline
$\xb_i$&{\small background - state vector }&$\yb_i$&output of observation operator ${\cal H}(\xb_i)$\\
\hline
$\bar x^b$&{\small average of $\xb_i$}&$\bar y^b$&{\small average of $\yb_i$}\\
\hline
$\Pb$&{\small background - error covariance}&&\\
\hline
$\xa$&{\small analysis - state vector}&$\Pa$& analysis - error covariance\\
\hline
\end{tabular}
\caption{Notations}
\label{table2}
\end{center}
\end{table}

\subsection{The UKF}

The unscented Kalman filter has been increasingly popular in engineering applications since its introduction about twenty years ago \cite{UKF:julier,julier2}. In a UKF algorithm, the error covariance is propagated with the dynamics using a set of vectors, or $\sigma$-points denoted by $\xsig$. Their definition is given in (\ref{eq_km1})-(\ref{eq_km2}). The $\sigma$-points are computed at each time step using a square root of the error covariance. In most UKF applications, $\sigma$-points are computed using either Cholesky factorization or matrix diagonalizations. In the notation, a variable with a superscript '$a$', such as $x^a$, represents the {\it analysis} value of the variable, i.e., the updated value based on observations. A variable with a superscript '$b$', such as $y^b$, represents the background, i.e.,  the propagated value of analysis using the system model. The algorithm is summarized as follows. At $t=k-1$, suppose we have the analysis and error covariance as well as its square root
\EQ
\label{eq_km1}
x^a (k-1), \; \Pa(k-1),\\
\Xa(k-1)=\sqrt{(n+\kappa)\Pa(k-1)},
\EE
where $\kappa$ is a scaling factor for the fine tuning of the higher order moments of the approximation error \cite{UKF:julier}. How to tune the value of $\kappa$ for a sparsity-based UKF is an open problem that needs further study. In this paper, $\kappa =0$ is used in all examples. A set of $\sigma$-points is generated as follows,
\EQ
\label{eq_km2}
\xsig_0(k-1)=x^a (k-1),\\
\xsig_i(k-1)=x^a(k-1)+\Xa_i(k-1), & 1\leq i\leq n, \\
\xsig_i(k-1)=x^a(k-1)-\Xa_i(k-1), & n+1\leq i\leq 2n. \\
\EE
The next step is to propagate the $\sigma$-points, which represent the background at $t=k$. For simplicity, the index '$k$' of all variables in the $k$th time-step is omitted.
\begin{equation}
\begin{aligned}
\xb_i&={\cal M}(\xsig_i(k-1)), & \yb_i&={\cal H}(\xb_i), &  0\leq i\leq 2n,\\
\bar x^b&= \displaystyle\sum_{i=0}^{2n} w_i\xb_i, & \bar y^b&=\displaystyle\sum_{i=0}^{2n} w_i\yb_i,
\end{aligned}
\end{equation}
where the weights are defined as follows
\EQ
w_0=\ds\frac{\kappa}{n+\kappa}, \; w_i=\frac{1}{2(n+\kappa)},
\EE
for $i=1,2,\cdots, 2n$. Define the variations
\EQ
\Xb_i=\xb_i-\bar x^b, \; \Yb_i=\yb_i-\bar y^b.
\EE
The background covariances are
\EQ
\Pb=\displaystyle \sum_{i=0}^{2n} w_i \Xb_i (\Xb_i)^T+Q,\\
P_{xy}=\displaystyle \sum_{i=0}^{2n} w_i \Xb_i (\Yb_i)^T,\\
P_{yy}=\displaystyle \sum_{i=0}^{2n} w_i \Yb_i (\Yb_i)^T+R.
\EE
The Kalman gain, $K$, satisfies the following equation,
\EQ
KP_{yy}=P_{xy}.
\EE
The analysis is updated as follows
\EQ
\begin{aligned}
\xa&=\bar x^b+K(y_o-\bar y^b),\\
\Pa &= \Pb -K(P_{xy})^T.
\end{aligned}
\EE
where $y_o$ is the observation at $t=k$. This completes one iteration of the filter. For the next step, $t=k+1$, go back to (\ref{eq_km1}) replacing the analysis by the updated value of $\xa$ and $\Pa$. 


\subsection{Sparse UKF}
The square root factorization of a matrix is not unique. For large and sparse matrices, various algorithms and their implementations on different computing platforms have been studied for many years. The literature can be traced back to the early days of electronic computers \cite{davis2}.  In the case of Cholesky factorization, the square root of a sparse matrix is still sparse, although the computation may require larger memory than the original matrix \cite{davis,rozin}. 

A dense error covariance is intractable in computation for global models used in NWP. In the proposed approach, we assume that $P$ and $\sqrt{P}$ are approximately sparse. In the algorithm, they are replaced by their sparse approximations, $\underline{P}$ and $\sqrt{\underline{P}}$. The indices of nonzero entries in the $i$th column are denoted by $\IP_i$ and $\ISP_i$, respectively. The sparsity index set of the forecast of $\sigma$-points, i.e. the background, are denoted by $\IBP$. Becasue  $\underline{P}$ and $\sqrt{\underline{P}}$ are approximations of $P$ and $\sqrt{P}$, the sparsity patterns do not have to be exact. In fact, in the example of  Lorenz-96 model presented in the next section, we assume that $\IP_i$, $\ISP_i$ and $\IBP$ equal to each other although $\sqrt{\underline{P}}$ may have a different sparsity pattern from that of $\underline{P}$.\\

\noindent \textbf{Algorithm I} (Sparse UKF)\\
Given the initial analysis,
\EQ
\xa(k-1), \; \spPa(k-1).
\EE

\noindent
\textbf{Step 1}. $\sigma$-points and forecast
\begin{equation}
\label{eq_sukf0}
\spXa(k-1)=\underline{\sqrt{(n+\kappa)\spPa(k-1)}},  \mbox{ sparsity index set }\ISP
\end{equation}
and
\begin{equation}
\label{eq_sukf1}
\begin{aligned}
\xb_0&={\cal M}(\xa(k-1)), & \yb_0&={\cal H}(\xb_0), \\
\spxb_{i}&={\cal M}(\xa(k-1)+\spXa_{i}(k-1); \IBP_i), & \yb_i&={\cal H}(\spxb_{i}\triangleright\xb_0), &&  1\leq i\leq n.\\
\spxb_i&={\cal M}(\xa(k-1)-\spXa_{i}(k-1); \IBP_i), & \yb_i&={\cal H}(\spxb_{i}\triangleright\xb_0), &&  n+1\leq i\leq 2n.\\
\bar x^b&=\displaystyle\sum_{i=0}^{2n} w_i(\spxb_i\triangleright\xb_0), &  \bar y^b&=\displaystyle\sum_{i=0}^{2n} w_i\yb_i\\
\end{aligned}
\end{equation}

\noindent \textbf{Step 2}. Background covariances
\EQ
\label{eq_sukf2}
\begin{aligned}
\spPb&=\displaystyle \sum_{i=0}^{2n} w_i\underline{(\spxb_i\triangleright\xb_0-\bar x)(\spxb_i\triangleright\xb_0-\bar x)^T}+\underline{Q}, &\mbox{ sparsity index set }{\cal I},\\
P_{xy}&=\displaystyle \sum_{i=0}^{2n} w_i(\spxb_i\triangleright\xb_0-\bar x)(\yb_i-\bar y)^T,\\
P_{yy}&=\displaystyle \sum_{i=0}^{2n} w_i(\yb_i-\bar y)(\yb_i-\bar y)^T+R.
\end{aligned}
\EE

\noindent \textbf{Step 3}. Kalman gain and analysis
\EQ
\label{eq_update}
\begin{aligned}
KP_{yy}&=P_{xy},\\
\xa&=\bar x^b+K(y_o-\bar y^b),\\
\spPa &= \spPb -\underline{K(P_{xy})}^T+\gamma I,& \mbox{ sparsity index set }{\cal I}.
\end{aligned}
\EE

In (\ref{eq_sukf2}), we assume $\spxb_0 = \xb_0$. The constant term $\gamma I$ in (\ref{eq_update}) is a diagonal matrix. The value of $\gamma$ is selected so that $\Pa$ is positive definite. The positive definiteness is guaranteed if $\gamma$ is larger than $|\lambda_{\min}|$, where $|\lambda_{\min}|$ is the smallest negative eigenvalue of 
\EQ
\label{eq_update2}
\spPb -\underline{K(P_{xy})}^T.
\EE
If the updated covariance matrix is positive definite, then $\gamma =0$. In all examples, the value of $\gamma$ is adaptively changed in every cycle depending on the smallest negative eigenvalue of (\ref{eq_update2})

\begin{figure}[!ht]
\begin{center}
\includegraphics[width=4.5in,height=2.5in]{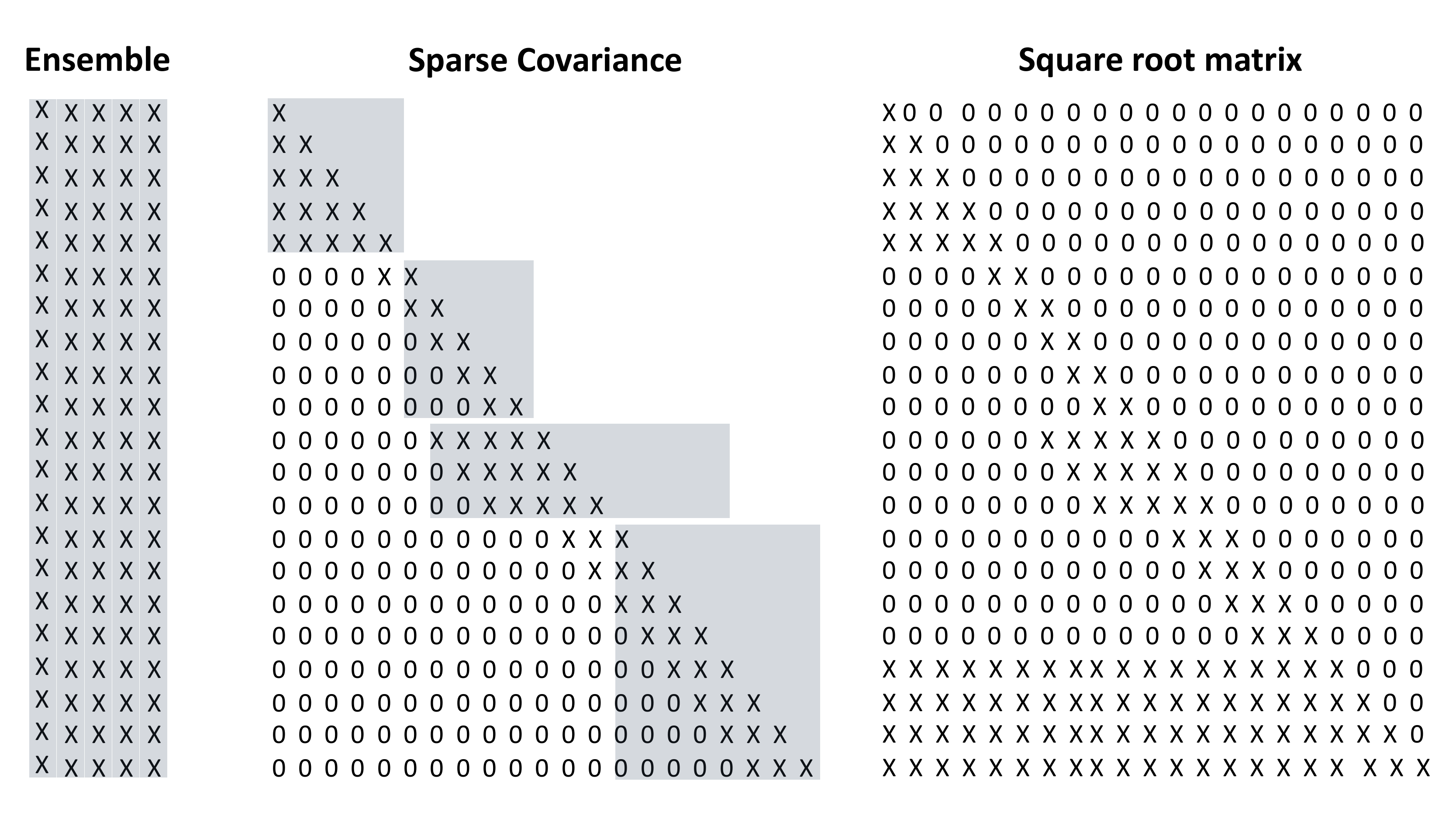} 
\caption{Patterns of ensemble vectors and sparse error covariances}
\label{fig_rank}
\end{center}
\end{figure}

The sparse UKF algorithm is based on the assumption that $\Pa$ can be  approximated by a sparse matrix $\spPa$. Although the $\sigma$-points in the algorithm play a similar role as that of ensembles in EnKF, using sparse UKF one can avoid the problem of rank deficiency. For systems with very high dimensions, the number of ensemble members used in an EnKF is much smaller than the dimension. As shown in Figure \ref{fig_rank}, the narrow and tall matrix of ensemble vectors makes EnKF fundamentally a rank deficient approach. In contrast, the block diagonal matrix $\spPa$ shown in Figure \ref{fig_rank} as a sparse approximation of $\Pa$ has full rank. If $N_{sp}$ of $\spPa$ is an integer close to the ensemble size of an EnKF, then the memory required by $\spPa$ is smaller than that by the ensemble matrix because of the symmetry of error covariance. The required computational load in Step 1 is extremely high if full state vectors are computed. Thanks to the sparsity, we only need to compute the entries with indices in $\IBP$. For a sparse UKF to be successful for problems with very high dimensions, it is critical to have component-based numerical models so that subsets of entries defined by $\IBP$ are computed in parallel; and most entries of state vectors are not evaluated at all. In addition each component-based computation requires a part of the state vector only. It saves the computer I/O load. 

Memory and (I/O) requirements are big factors influencing the efficiency of computational algorithms. Because an error covariance and its diagonal blocks are symmetric, the memory and I/O usage can be reduced by almost a half for subproblems with symmetry. The granularity of a computational algorithm has considerable impact on its efficiency in parallel computing. By granularity control we mean that one can divided a high dimensional problem into subproblems of desired dimensions. As shown in Figure \ref{fig_rank}, the error covariance and its square root consist of sparse blocks or sparse columns. This is different from the EnKF in which state vectors in an ensemble are dense. In a sparse approximation of $\Pa$, the number of nonzero entries, a parameter similar to that used for distance-based localization methods, can be easily changed in $\underline{\Pa}$ so that the error covariance and its square root can be grouped into smaller blocks of different sizes for parallel computation. 

\subsection{Lorenz-96 model}
\begin{figure}[!ht]
\begin{center}
\includegraphics[width=3.0in,height=2.5in]{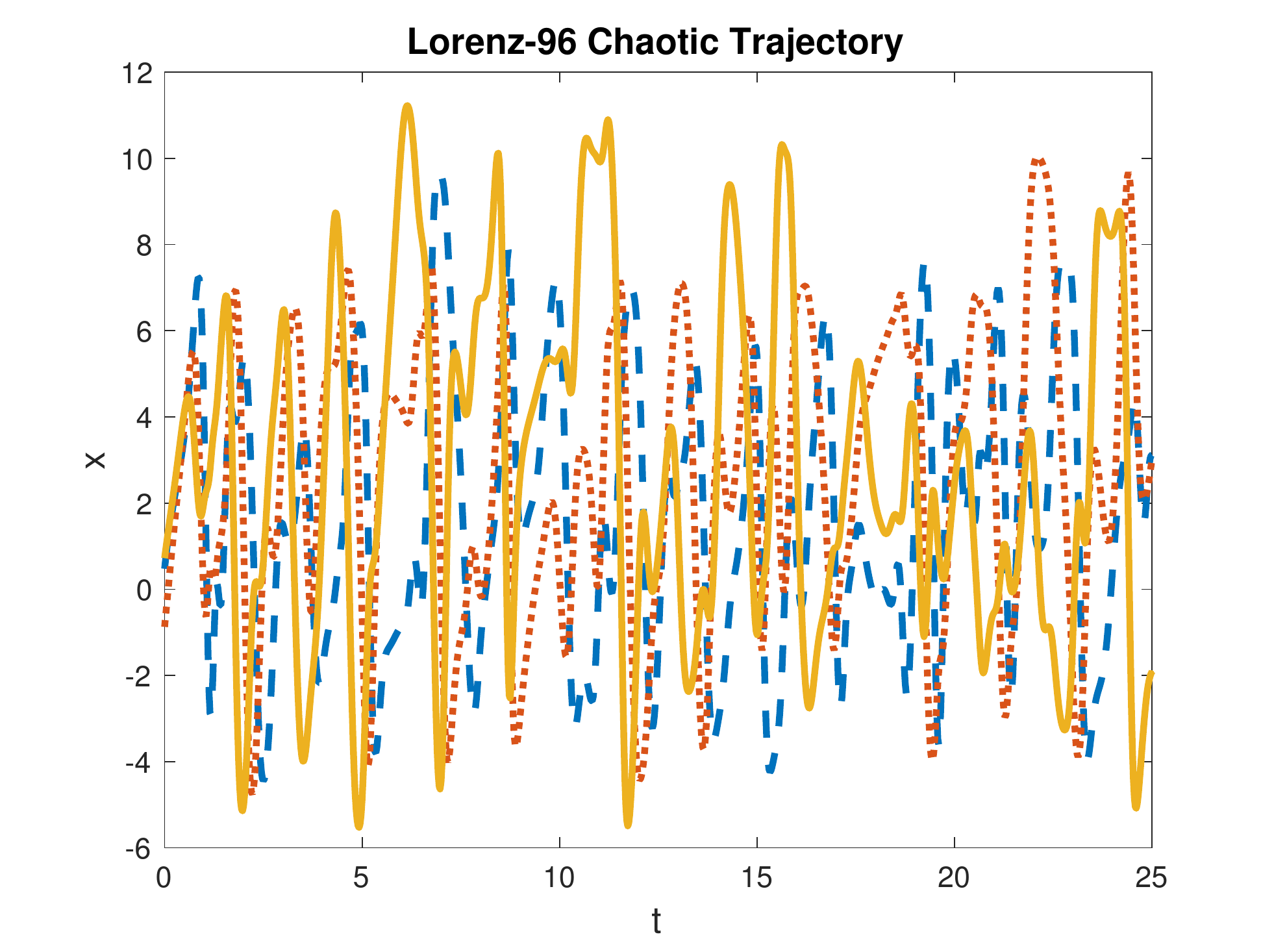} 
\caption{A chaotic trajectory of the Lorenz-96 model, $x_1$(solid), $x_2$(dash), $x_3$(dot).}
\label{fig_chaos}
\end{center}
\end{figure}

In this section, we use a Lorenz-96 model that was first introduced in \cite{lorenz96} to test the performance of the sparse UKF. Consder
\EQ
\label{lorenz_1}
\begin{aligned}
\Fr{dx_i}{dt}&=(x_{i+1}-x_{i-2})x_{i-1}-x_i+F, & i=1, 2, \cdots, n,\\
x_{n+1}&=x_1,\\
n&=40,\\
\Delta t&=0.025,\\
 F&= 8.
 \end{aligned}
\EE
The system has chaotic trajectories as shown in Figure \ref{fig_chaos}, a plot of $x_1(t), x_2(t), x_3(t)$.  
The simulations are conducted based on a 4th-order Runge-Kutta discretization. The trajectories are used as the ground truth. The sparsity pattern for $\spPa$ and $\sqrt{\spPa}$ are assumed to be centered along the diagonal line with a fix length. The total number of nonzero entries in each column is $N_{sp}$. We would like to point out that the sparse matrices are approximations of the true error covariance and its square root. The sparsity pattern of $\sqrt{\spPa}$ is, in fact, different from that of $\spPa$. In the approximation, however, we ignore the difference and use the same sparsity pattern for both. As a result, the memory required by $\sqrt{\spPa}$ is reduced almost by half. This idea works fine for the Lorenz-96 model. However, a systematic way of choosing the sparsity pattern for $\sqrt{\spPa}$ based on given $\spPa$ is an open problem that needs further study. 

\begin{figure}[!ht]
\begin{center}
\includegraphics[width=3.0in,height=2.0in]{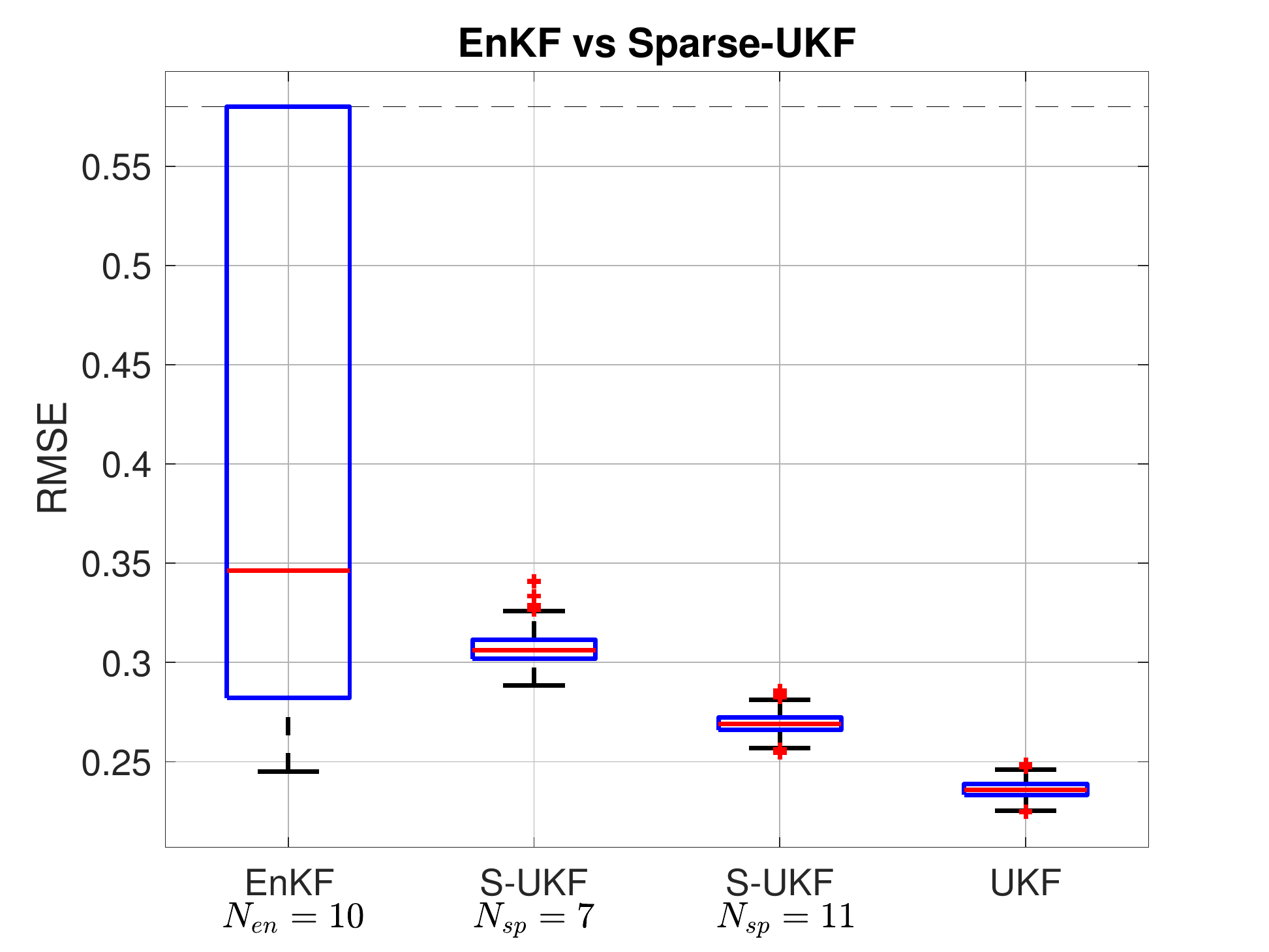} 
\caption{Boxplot of RMSE}
\label{fig_box_UKF}
\end{center}
\end{figure}

The numerical experimentation is based on $N=1000$ uniformly distributed random initial states in $[-1\; 1]$. The time step size is $\Delta t=0.025$. The total number of time steps for each simulation is $N_t=4000$.  The number of observations at any given time is $m=20$, i.e. every other state variable is measured. The observational error has the Gaussian distribution with a covariance $R=I$, the identity matrix. The initial background error covariance is $\Pb(0)=0.2I$. The following RSME is used to measure the accuracy of estimation
\EQ
RSME=\sqrt{\Fr{1}{n(N_t+1)}\displaystyle\sum_{k=0}^{N_t}||\xa(k)-x^{truth}(k)||^2_2}.
\EE
For comparison, an EnKF is applied to the same data set. The localization radius is $\rho=4$ and the inflation factor is $\sqrt{1.08}$. A full scale UKF based on dense error covariance is also applied as the best estimator for comparison. As an indication of computational load, the number of entries to be computed in each algorithm is shown in Table \ref{table_1}. The boxplot of simulation results is shown in Figure \ref{fig_box_UKF}. To summarize, the sparse UKF has considerably smaller error variation than that of the EnKF. This is expected because the new approach avoids the problem of rank deficiency. The medians of estimation errors are also smaller than that of EnKF. In the cases of $N_{sp}=7$ and $11$, the memory size required by the sparse error covariance is smaller than the memory size needed to store the ensemble vectors in the EnKF if  $N_{ens}=10$. The Cholesky factorization, while maintaining the sparsity property, may require additional memory. In this example of Lorenz-96, we use the same sparsity pattern for both $\sqrt{\spPa}$ and $\spPa$. This assumption simplifies the algorithm and reduces the memory and I/O requirements. The computational load, in the number of entry evaluations, is increased for sparse UKF because of the number of $\sigma$-points is $2n$. Reducing the number of entries being evaluated and testing the impact of Cholesky factorization on the efficiency of UKF  is ongoing research that is not addressed in this paper.  

\begin{table}[h]
\begin{center}
\begin{tabular}{|c|c|c|c|c|c|}
\hline
{\small Filter}& {\small Size}& {\small Entries}&{\small Error} &{\small Error}&{\small Error}\\
&&{\small  EVAL}&{\small Median}&{\small Mean}&{\small STD}\\
\hline
{\small EnKF}&{\small $N_{ens}=10$}& {\small 400}&0.3462&1.0741&1.0652\\
\hline
{\small S-UKF}&{\small $N_{sp}=7$}&{\small 600}&0.3061&0.3067&0.0071\\
\hline
{\small S-UKF}&{\small $N_{sp}=11$}&{\small 920}&0.2691&0.2691&0.0048\\
\hline
\end{tabular}
\caption{Summary of simulation results}
\label{table_1}
\end{center}
\end{table}

\section{Progressive EKF}

In a sparse UKF, the $\sigma$-points are computed by taking a square root of the error covariance, such as the Cholesky factorization. The resulting $\sigma$-points are sparse. However, this process may require additional memory and computation. In this section, we propose a progressive algorithm of approximating error covariance without taking square roots. 

\subsection{Basic ideas}
The main assumption for this algorithm is the following progressive relationship  
\EQ
\label{eq_prog}
M_{k-1}\Pa(k-1) M_{k-1}^T=\Pa (k-1) + \Delta \Pb,
\EE
where $\Delta \Pb$ is assumed to be small. In (\ref{eq_prog}), $M_{k-1}$ is the Jacobian of ${\cal M}$ at $\xa (k-1)$. Similarly, the Jacobian of ${\cal H}$ is $H_k$. 
To estimate $\Delta \Pb$, assume
\EQ
\label{eq_DM}
M_{k-1}&=I+\Delta M_{k-1}.
\EE
where we assume that $\Delta M_{k-1}$ is small. If the system model is based on the discretization of a differential equation with a small time step size, then
\EQ
\label{eq_prog3}
{\cal M}(x(k-1))=x(k-1) + O(\Delta t^\alpha), \;\; \alpha \geq 1.
\EE
The Jacobian of $O_{k-1}(\Delta t^\alpha)$ in space variables is expected to have small value if $\Delta t$ is small, which makes (\ref{eq_DM}) a reasonable assumption. Then we have 
\EQ
\label{eq_prog2a}
\begin{aligned}
&M_{k-1}\Pa(k-1) M_{k-1}^T\\
&= (I+\Delta M_{k-1})\Pa(k-1)(I+\Delta M_{k-1}^T)\\
&= \Pa(k-1)+\Delta M_{k-1} \Pa(k-1)+\left( \Delta M_{k-1} \Pa(k-1)\right)^T+ \Delta M_{k-1} \Pa(k-1)\Delta M_{k-1}^T\\
&\approx \Pa(k-1)+\Delta M_{k-1} \Pa(k-1)+\left( \Delta M_{k-1} \Pa(k-1)\right)^T.
\end{aligned}
\EE
This is in consistent with (\ref{eq_prog}). It can be computed using a tangent linear model. Or it can be approximated using the dynamical model
 \EQ
\label{eq_prog2b}
\begin{aligned}
&M_{k-1}\Pa(k-1) M_{k-1}^T\\
&= (I+\Delta M_{k-1})\Pa(k-1)(I+\Delta M_{k-1}^T)\\
&\approx \left({\cal M}(x(k-1)+\delta \Pa(k-1))-{\cal M}(x(k-1))\right)/\delta\\
&+\left({\cal M}(x(k-1)+\delta \Pa(k-1))-{\cal M}(x(k-1))\right)^T/\delta-\Pa.
\end{aligned}
\EE
where $\delta>0$ is the step size of a finite difference approximation of $\Delta M_{k-1}\Pa$. Its value should be determined depending on the numerical model and its linearization. In (\ref{eq_prog2b}), a vector and matrix summation is a new matrix resulting from adding the vector to every column in the matrix. Applying an operator to a matrix is to apply the operator to every column in the matrix. 

\subsection{Progressive EKF}
The column vectors in the matrices in (\ref{eq_prog2a}) and (\ref{eq_prog2b}) are sparse. However, the number of column vectors equals $n$, which can be as high as $10^6-10^7$ for some atmospheric models. Applying a full model to all the vectors is impractical  because of the high computational and I/O loads.  On the other hand, if we approximate the error covariance using a given sparsity, only a small portion of the entries in each column vector is to be evaluated. Evaluating the entire state vector is unnecessary. This is the reason we need a component-based model. Then the algorithm of progressive EKF is summarized as follows. \\

\noindent \textbf{Algorithm II} (Progressive EKF)\\
Given the initial analysis at $t=k-1$,
\EQ
\xa(k-1) \mbox{ and }\spPa(k-1).
\EE

\noindent
\textbf{Step 1}. Forecast
\EQ
\xb={\cal M}(\xa(k-1)),\\
\yb={\cal H}(\xb).
\EE

\noindent
\textbf{Step 2}. Background error covariance
\EQ
\label{eq_prop4}
\begin{aligned}
\spPb&=\underline{\left({\cal M}\left(\xa(k-1)+\delta \spPa(k-1),{\cal I}\right)-\xb \right)}/\delta\\
&+\underline{\left({\cal M}\left(\xa(k-1)+\delta \spPa(k-1),{\cal I}\right)-\xb \right)}^T/\delta-\spPa+Q.
\end{aligned}
\EE

\noindent
\textbf{Step 3}. Kalman gain and analysis
\EQ
\begin{aligned}
K&=\spPb H_{k}^T(H_{k}\spPb H^T_{k} +R)^{-1},\\
\xa&= \xb+K(y_o- \yb),\\
\spPa&= (I-KH_{k})\spPb.\\
\end{aligned}
 \EE

Different from the sparse UKF, this algorithm avoids the computation of matrix square roots. However, the algorithm requires that $\Delta \Pb$ in (\ref{eq_prog}) can be approximated effectively. From (\ref{eq_prog3}), the method is expected to work better for a small time step-size. If $\Delta t$ is large, $\Delta M_{k-1}$ in (\ref{eq_DM}) may not be small enough. A remedy is to use a refined step-size in an inner-loop computation. More specifically, the discrete model is a discretization of a continuous-time model. The discrete time moment $k-1$ corresponds to the continuous time moment $(k-1)\Delta t$. 
We refine the step size by dividing the time interval into $n_p$ subintervals. In our examples, we choose $n_p=2$. The refined time steps are 
\EQ
(k-1)\Delta t, (k-1)\Delta t+\Fr{\Delta t}{n_p}, \cdots, (k-1)\Delta t+s\Fr{\Delta t}{n_p},\cdots, k\Delta t,  \;\;\;\;\; 0\leq s \leq n_p
\EE
For the inner loop,  one can compute a sequence of backgrounds, $\tilde x^b (s)$. 
\EQ
\label{refinestep}
t_s=(k-1)\Delta t+s\Fr{\Delta t}{n_p},\\
\tilde x^b(s)=\tilde{\cal M}_{t_s} (\xa(k-1)), & s=1,2, \cdots, n_p.
\EE
where $\tilde {\cal M}_{t_s}$ represents the refined model function in the time interval from $t=(k-1)\Delta t$ to $t=t_s$. In Step 2, repeat (\ref{eq_prop4}) $n_p$ times along the sequence of background states, $\tilde x^b(s)$, without adding $Q$ until the last round. This refined Step 2 increases the computational load, while improving the accuracy of the progressive estimation. 

\subsection{Examples}
\begin{figure}[!ht]
\begin{center}
\includegraphics[width=3.0in,height=2.0in]{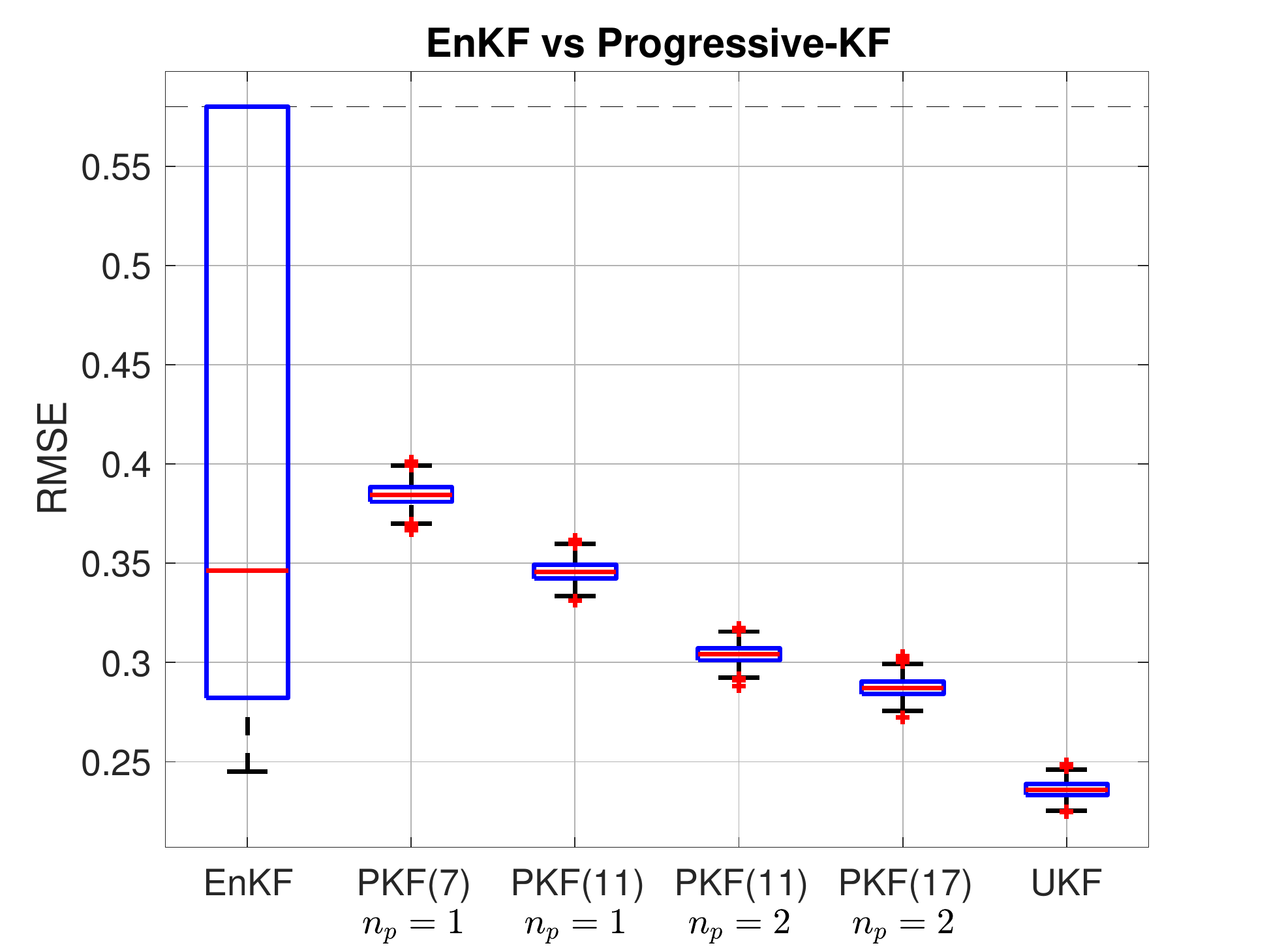} 
\caption{Boxplot of RMSE. For Prograssive-KF, $N_{sp}$ = 7, 11, and 17}
\label{fig_progEKF}
\end{center}
\end{figure}

In the following, we apply the progressive EKF to the Lorenz-96 model using the same parameters as in (\ref{lorenz_1}). The error covariance is approximated using sparsity matrices with $N_{sp}=7, 11, 17$. For $N_{sp}=11$, we tested the idea of refining step-size using $n_p=1$ and $n_p=2$. The results are shown in Figure \ref{fig_progEKF} and summarized in Table \ref{table2}.  Comparing to EnKF, the error variations of the progressive EKF are uniformly and significantly smaller. If $N_{sp}=7$, which is smaller than the ensemble size $N_{ens}=10$, the median value of estimation error is larger than that of the EnKF. The median error for $N_{sp}=11$ is comparable to that of the EnKF. If a refined step-size in (\ref{refinestep}) is applied, for instance $n_p=2$, the median estimation error is further reduced. Comparing to the performance of the sparse UKF, the error variations are similar. However, the estimation error of the sparse UKF has a smaller median in all cases. For example, to achieve a similar performance as the sparse UKF when $N_{sp}=11$,  one has to use a larger sparsity index $N_{sp}=17$ for the progressive EKF. 
The computational load of the progressive EKF, in terms of the component evaluation, is comparable to the EnKF, as shown in Table \ref{table2}. 
\begin{table}[h]
\begin{center}
\begin{tabular}{|c|c|c|c|c|c|}
\hline
{\small Filter}& {\small Size}& {\small Entries}&{\small Error} &{\small Error}&{\small Error}\\
&&{\small  EVAL}&{\small Median}&{\small Mean}&{\small STD}\\
\hline
{\small EnKF}&{\small $N_{ens}=10$}& {\small 400}&0.3462&1.0741&1.0652\\
\hline
{\small P-KF}&{\small $N_{sp}=7$}&&&&\\
&{\small $N_p=1$}&{\small 320}&0.3845&0.3846&0.0055\\
\hline
{\small P-KF}&{\small $N_{sp}=11$}&&&&\\
&{\small $N_p=1$}&{\small 480}&0.3455&0.3458&0.0050\\
\hline
{\small P-KF}&{\small $N_{sp}=11$}&&&&\\
&{\small $N_p=2$}&{\small 480x2}&0.3041&0.3041&0.0044\\
\hline
{\small P-KF}&{\small $N_{sp}=17$}&&&&\\
&{\small $N_p=3$}&{\small 720x2}&0.2872&0.2873&0.0046\\
\hline
\end{tabular}
\caption{Summary of simulation results}
\label{table2}
\end{center}
\end{table}

\section{Conclusions}
Two algorithms of Kalman Filters based on sparse error covariances are introduced. They are tested using the Lorenz-96 model with $40$ state variables and chaotic trajectories. Both algorithms share the same basic idea: the error covariance is approximated using a sparse matrix. Thanks to the sparsity, the required memory size is significantly reduced. The symmetry of the error covariance can potentially reduce the I/O load. The analysis error covariance can be updated as a sparse matrix in each cycle using a deterministic process, either a square root matrix or a progressive algorithm. The updated sparse matrix is then used as the background error covariance for the next cycle. Relative to the EnKF, the main advantage of the proposed methods is that the estimation process do not need an ensemble; and the error covariance has a full rank. The algorithms do not suffer issues of rank deficiency as in EnKFs. As a result, the variation of analysis error is constantly small in all examples. Techniques of localization and covariance inflation are unnecessary.  Relative to 4D-Var methods, the proposed algorithms are mostly parallel. They provide not only the state estimate but also the analysis error covariance. For the purpose of scalability, we suggest that the proposed methods are applied with component-based numerical models. From the examples, the sparse UKF has better accuracy than the progressive EKF. If the computational load of taking square roots of sparse matrices is affordable, then the sparse UKF is the approach of our choice. On the other hand, the progressive EKF is a simple algorithm that avoids taking square roots of large matrices, provided that the progressive approximation of error covariance is adequately accurate. Although most conclusions drawn in this paper are based on simulations using the Lorenz-96 model, the algorithms are developed for general applications. Testing the methods using different types of system models is a main topic of our future work.


\begin{thebibliography}{99}
\bibitem{davis} T. A. Davis, Direct Methods for Sparse Linear Systems, SIAM, 2006.


\bibitem{davis2} T. Davis, S. Rajamanickam, and W. M. Sid-Lakhdar, A survey of direct methods for sparse linear systems, Acta Numerica, Vol. 25, 2016, pp. 383-566.


\bibitem{houtekamer} P. L. Houtekamer and F. Zhang, Review of the Ensemble Kalman Filter for Atmospheric Data Assimilation, Monthly Weather Review, Vol. 144, 2016, pp. 4489 - 4532.

\bibitem{UKF:julier} S. Julier, J. Uhlmann, and H. F. Durrant-Whyte, A New Method for the Nonlinear Transformation of Means
and Covariances in Filters and Estimators, {\it IEEE Trans. Automatic Control}, Vol. 45 (3), 2000, pp. 477-482.

\bibitem{julier2} S. J. Julier and J. K. Uhlmann, Unscented filtering and nonlinear
estimation, {\it Proceedings of the IEEE,} Vol. 92, No. 3, 2004, pp.

\bibitem{lorenz96} E. Lorenz, Predictability -- A problem partly solved. Seminar on Predictability, Vol. I, ECMWF, 1996.

\bibitem{rozin} E. Rozin and S. Toledo, Locality of reference in sparse Cholesky factorization methods, Electronic Transactions on Numerical Analysis, Vol. 21, 2005, pp. 81-106.

\bibitem{xu} L. Xu, R. Rosmond, and R. Daley, Development of NAVDAS-AR: Formulation and initial tests of the linear problem, Tellus, 57A, 2005, pp. 546-559.
\end{thebibliography}
\end{document}